\documentclass[10pt,reqno]{amsart}
\usepackage{amsfonts}
\usepackage{amsmath}
\usepackage{amssymb}
\usepackage{color}
\usepackage{epsfig,exscale}
\usepackage[all]{xy}
\usepackage{graphicx}
\usepackage{xspace}

\newtheorem{prop}{Proposition}[section]
\newtheorem{lem}{Lemma} [section]
\newtheorem{cor}{Corollary}[section]

\newtheorem{thm}{Theorem}[section]
\newtheorem{de}{Definition}[section]

\DeclareMathOperator{\End}{End}       
\DeclareMathOperator{\Max}{Max}       
\DeclareMathOperator{\Prim}{Prim}     
\DeclareMathOperator{\Sup}{Sup}       

\newcommand{\A}{\mathcal{A}}        
\newcommand{\dl}{\delta}            
\newcommand{\id}{\mathrm{id}}       
\newcommand{\K}{\mathbb{K}}         
\newcommand{\N}{\mathbb{N}}         
\newcommand{\Rt}{\tilde{R}}         
\renewcommand{\th}{\theta}          

\def\<#1,#2>{\langle#1,#2\rangle}   
\newcommand{\sepword}[1]{\quad\mbox{#1}\quad} 

\begin{document}

\title{Rota--Baxter algebras and new combinatorial identities}

\author{Kurusch Ebrahimi-Fard}
\address{Max-Planck-Institut f\"ur Mathematik,
         Vivatsgasse 7,
         D-53111 Bonn, Germany.}
 \email{kurusch@mpim-bonn.mpg.de}
 \urladdr{http://www.th.physik.uni-bonn.de/th/People/fard/}

\author{Jos\'e M. Gracia-Bond\'{\i}a}
\address{Departamento de F\'{\i}sica Te\'orica I,
         Universidad Complutense,
         Madrid 28040, Spain}

\author{Fr\'ed\'eric Patras}
\address{Laboratoire J.-A. Dieudonn\'e
         UMR 6621, CNRS,
         Parc Valrose,
         06108 Nice Cedex 02, France}
 \email{patras@math.unice.fr}
 \urladdr{www-math.unice.fr/~patras}

\date{January 10, 2007}


\begin{abstract}
The word problem for an arbitrary associative Rota--Baxter algebra
is solved. This leads to a noncommutative generalization of the
classical Spitzer identities. Links to other combinatorial
aspects, particularly of interest in physics, are indicated.
\end{abstract}

\maketitle


\keywords


{\footnotesize{\noindent Keywords: Rota--Baxter relation; free
algebras; word problem; quasi-symmetric functions; noncommutative
symmetric functions; Hopf algebra; pre-Lie relation; Dynkin
idempotent; Spitzer's identity; Bohnenblust--Spitzer identities;
Magnus' expansion}}

\smallskip

{\footnotesize{\noindent\textsl{Mathematics Subject Classification
2000}}: 05A19,16W30}


\thispagestyle{empty}


\section{Introduction and definitions}
\label{sect:Introibo}

Nearly forty years ago, a class of combinatorial formulas for
random variables were recast by Rota as identities in the theory
of Baxter maps~\cite{baxter1960}. The key result was the solution
of the word problem, for associative, commutative algebras endowed
with such maps. This showed the equivalence of the combinatorics
of fluctuations with that of classical symmetric functions. Since
then, operators of the Baxter type kept showing up in all sorts of
applications, and lately in the Hopf algebraic approach to
renormalization~\cite{eg2006}.  In many instances, the algebra in
question is not commutative.  The time has come to revisit the
word problem, and the corresponding identities, in the
noncommutative case. Roughly speaking, we are led to replace
symmetric functions of commuting variables by quasi-symmetric
functions of non-commuting ones. Sequences of `noncommutative
Spitzer identities' ensue. In an applied vein, we explore the
connection of our word problem with Lam's approach to the Magnus
expansion for ordinary differential equations.

\begin{de}
Let $\K$ be a field of characteristic zero.  Let $A$ be a
$\K$-algebra, not necessarily associative nor commutative nor
unital. An operator $R \in \End(A)$ satisfying the relation
\begin{equation}
Ra \, Rb = R(Ra\,b + aRb) + \theta R(ab), \sepword{\rm for all} a,b
\in A, \label{eq:RBR}
\end{equation}
is said Rota--Baxter of weight~$\th \in \K$. The pair~$(A,R)$ is a
weight~$\th$ \textit{Rota--Baxter algebra} {\rm(RBA)}.
\end{de}
The Rota--Baxter identity~\eqref{eq:RBR} prompts the definition of a
new product $a \ast_R b :=Ra\,b + aRb + \theta ab$, $a,b \in A$.

\begin{prop}
The linear space underlying~$A$ equipped with the product~$\ast_R$
is again a {\rm RBA} of the same weight with the same Rota--Baxter
map. We denote it by $(A_R,R)$.  If~$A$ is associative, so
is~$A_R$.
\end{prop}

We call $\ast_R$ the Rota--Baxter double product. Clearly $R$
becomes an algebra map from~$A_R$ to~$A$. Note that $\Rt :=-\theta
\id_A - R$ is Rota--Baxter as well, and $*_{\tilde R}=-*_R$. One
may think of Rota--Baxter operators as generalized integrals.
Indeed, relation~\eqref{eq:RBR} for the weight $\th=0$ corresponds
to the integration-by-parts identity for the Riemann integral; the
reader will have no difficulty in checking duality
of~\eqref{eq:RBR} with the `skewderivation' rule
\begin{equation*}
\dl(ab) = \dl a\,b + a\dl b + \th\dl a\dl b.
\end{equation*}
For instance, the finite difference operator of step~$-\theta$,
given by $\dl f(x):=\theta^{-1}(f(x-\theta)-f(x))$, is a
skewderivation. The summation operator $Zf(x):=\sum_{n\geq
1}\theta f(x + \theta n)$ is Rota--Baxter of weight~$\theta$, and
we find $\dl Z=\id=Z\dl$ on suitable classes of functions. Scaling
$R \to \th^{-1}R$ reduces the study of {\rm RBA}s of nonvanishing
weight to the case $\th=1$.  For notational simplicity we proceed
considering this one, returning to general weight when convenient.
Also, henceforth we assume we are dealing with associative {\rm
RBA}s; non-associative {\rm RBA}s will arise later in an ancillary
role.


\section{Main result}
\label{sect:Rota-bis}

We now extend to our noncommutative setting Rota's notion of
\textit{standard} RBA, see~\cite{cartier1972,rota69,rotasmith72}.
Let $X=(x_1,\dots,x_n,\dots)$ be a countably infinite, ordered set
of variables and $T(X)$ the tensor algebra over~$X$. The elements
of $X$ are called noncommutative polynomials (over $X$). Consider
the pair $(\A,\rho)$, where $\A$ is the algebra of countable
sequences $\Upsilon\equiv(y_1,\dots, y_n,\dots)$ of elements
$y_i\in T(X)$ with pointwise addition and product, and $\rho$
given by \begin{equation*} \rho\Upsilon = (0, y_1, y_1 + y_2, y_1
+ y_2 + y_3, \dots).  \end{equation*} By abuse of notation we
regard $X$ itself as an element of~$\A$.  The component~$y_p$
of~$\Upsilon$ is denoted~$\Upsilon_p$.

\begin{lem}
The algebra $\A$ together with $\rho \in \End(\A)$ defines a
weight $\theta=1$ Rota--Baxter algebra structure.
\end{lem}

This is a straightforward verification.  We remark that~$\rho$ has a left
inverse.

\begin{thm}
The Rota--Baxter subalgebra $(\mathcal R,\rho)$ of~$\A$ generated by~$X$ is
free on one generator in the category of $\K$-{\rm RBAs}.
\end{thm}

In detail, our assertions are the following.

\begin{itemize}
\item{} $X \in \mathcal R$.

\item{} The product in~$\mathcal R$ is associative.

\item{} $\rho$ is a Rota--Baxter operator.

\item{} Let $(A,R)$ be any associative RBA and $a\in A$. There is
a unique algebra map $h:\mathcal R\to A$ with $h(X)=a$ and such
that $R \circ h=h \circ \rho$.
\end{itemize}

The pair $(\mathcal R,\rho)$ is what we call the standard RBA. The
point of course is that the theorem allows us to prove the validity
for any RBA $A$ of an identity involving one element of~$A$ and~$R$,
by proving it for~$X$ in~$\mathcal R$.

Only the last assertion in the list above asks for proof.  We shall
follow Rota and Smith~\cite{rotasmith72} as far as possible.  The
adaptation to the noncommutative setting requires a bit of care.  The
lexicographical ordering $<_L$ for noncommutative monomials over~$X$
is useful; for any noncommutative polynomial~$P$ we write $\Sup P$ for
the highest monomial in~$P$ for~$<_L$ and extend the lexicographical
ordering of noncommutative monomials to a partial ordering on $T(X)$.
Namely, we write $P<_LP'$ whenever $\Sup P<_L\Sup P'$.  Note that, for
$P,P'$ homogeneous noncommutative polynomials and $z,t$ in~$T(X)$, we
have
$$
P <_L P' \Rightarrow Pz<_L P'z \sepword{\rm and} z <_L t \Rightarrow
Pz <_L Pt.
$$
Henceforth we just employ the generic $R$ for the Rota--Baxter map on
the standard RBA; this should not lead to any confusion.

\begin{proof}
(Main steps.)  Let us call $\End$-algebra any associative algebra
$W$ provided with a distinguished endomorphism $T_W$, so that an
$\End$-algebra morphism $f$ from~$W$ to~$W'$ satisfies $f \circ
T_W=T_{W'}\circ f$.  Write~$\mathcal L$ for the free
$\End$-algebra on one generator~$Z$. The elements of $\mathcal L$
are linear combinations of all symbols obtained from~$Z$ by
iterative applications of the endomorphism~$T$ and of the
associative product; they look like $ZT^2(TZ \, T^3Z)$, and so on.
We call these symbols $\mathcal L$-monomials. A RBA~$A$ is an
$\End$-algebra together with the relation~\eqref{eq:RBR} on $T_A
\equiv R$.  Denote by~$\mathcal F$ the free RBA on one
generator~$Y$. Between the three algebras $\mathcal L$, $\mathcal
F$, $\mathcal R$ there are the following maps: unique
$\End$-algebra maps $F$, $U$ from~$\mathcal L$ to~$\mathcal F$,
respectively~$\mathcal R$, sending~$Z$ to~$Y$ respectively~$X$;
and a unique onto Rota--Baxter map $h'$ sending~$Y$ to~$X$.
Moreover $U=h'\circ F$.

We have to show the existence of an inverse for~$h'$ in the RBA
category.  Clearly $\ker F\subseteq\ker U$.  We need only prove that
$\ker U\subseteq\ker F$.

Any $l\in\mathcal L$ can be written uniquely as a linear combination
of $\mathcal L$-monomials.  We write $\Max l$ for the maximal number
of $T$'s occurring in the monomials, so that, say, $\Max(ZT^2(ZTZ)+
Z^3T^2Z\,Z)=3$.  We call~$\alpha$, a $\mathcal L$-monomial,
\textit{elementary} iff it can be written as either $Z^i,i\geq 0$ or
as a product $Z^{i_1}Tb_1 \, Z^{i_2} \cdots Tb_k\, Z^{i_{k+1}}$, where
the $b_i$s are elementary, and $i_2,\ldots,i_k$ are strictly positive
integers, while $i_1$ and $i_{k+1}$ may be equal to zero; this
definition makes sense by induction on $\Max\alpha$.  It turns out
that every element~$l$ of~$\mathcal L$ can be written as the sum of a
linear combination of elementary monomials with an element $r_l$ such
that $F(r_l)=0$.  This is due to the fact that, up to the addition of
suitable elements in $\ker F$, products like $Tc \, Td$ can be
iteratively cancelled from the expression of~$l$ using
relation~\eqref{eq:RBR}.

We claim that for $p$ large enough and $l\not=l'$, with $l,l'$
elementary monomials, we have $\Sup U(l)_p\not=\Sup U(l')_p$, from
which the required $\ker U\subseteq\ker F$ follows.  Our assertion can
be verified by induction on $\Max l$, using that~$U$ is an
$\End$-algebra map.
\end{proof}

\begin{cor}
The images of the elementary monomials of~$\mathcal L$ in~$\mathcal R$
form a linear basis of the free RBA on one generator.
\end{cor}


\section{Two interesting Hopf algebras}
\label{sect:resume}

Inductively define in a general RBA $(A,R)$,
\begin{equation*}
(Ra)^{[n+1]} = R\big((Ra)^{[n]}a\big) \sepword{\rm and} (Ra)^{\{n+1\}} =
R\big(a(Ra)^{\{n\}}\big).
\end{equation*}
with the convention that $(Ra)^{[1]}=Ra=(Ra)^{\{1\}}$ and
$(Ra)^{[0]}=1=(Ra)^{\{0\}}$, with the unit adjoined if need be.
These iterated compositions with~$R$ appear in the context of
Spitzer formulas.  Of course there is no difference between
$(Ra)^{[n]}$ and $(Ra)^{\{n\}}$ in the commutative context.

Coming back to the standard {\rm{RBA}} $(\mathcal R,R)$, notice that:
\begin{equation*}
\big(R(y_1, y_2, y_3,\dots)\big)^{[2]} = R\bigl(R(y_1, y_2,
y_3,\dots)\, (y_1, y_2, y_3,\dots)\bigr)= (0, 0, y_1y_2, y_1y_2 +
y_1y_3 + y_2y_3,\dots)
\end{equation*}
This begins to give the game away.  In general, the $(n+1)$-th
entry of $\big(R(y_1, y_2, y_3,\dots)\big)^{[k]}$ is the
elementary `symmetric' function of degree~$k$, restricted to the
first~$n$ variables, the $(n+2)$-th entry is given by the same,
restricted to~$n+1$ variables, and so on.  The quotes on
`symmetric' remind us here that the $y_i$ do not commute. The
pertinent notion here is Hivert's quasi-symmetric functions over a
set of noncommuting variables~\cite{bergeron05,novelli06}.  Denote
as usual by~$[n]$ the set of integers between~$1$ and~$n$. Let $f$
be a surjective map from~$[n]$ to~$[k]$.  Then the quasi-symmetric
function $M_f$ over~$X$ associated to~$f$ is by definition
\begin{equation*}
M_f\,X = \sum\limits_\phi x_{\phi^{-1}\circ f(1)} \cdots
x_{\phi^{-1}\circ f(n)},
\end{equation*}
where $\phi$ runs over the set of increasing bijections between
subsets of $\N$ of cardinality $k$ and $[k]$.  Let us
represent~$f$ as the sequence of its values, $f=f(1),\dots,f(n)$,
in the notation $M_f$. We also denote by $M_f^{l}$ the image
of~$M_f$ under the map sending $x_i$ to~$0$ for $i>l$ and to
itself otherwise. For example,
$$
M_{1,3,3,2}\,X = x_1x_3x_3x_2 + x_1x_4x_4x_2 + x_1x_4x_4x_3 +
x_2x_4x_4x_3 + \dots \sepword{\rm and} M^3_{1,3,3,2}\,X =
x_1x_3x_3x_2.
$$
The linear span ${\rm NCQSym}(X)$ of the $M_f$ ---a subalgebra of
the completion of the algebra of noncommutative polynomials
over~$X$--- is related to the Coxeter complex of type $A_n$ and
the corresponding Solomon--Tits and twisted descent
algebras~\cite{ps2005}.

Finally, write $[n]$ for the identity map on~$[n]$ and~$\omega_n$
for the endofunction of~$[n]$ reversing the ordering, so that
$M_{\omega_n}=M_{n,n-1,\dots ,1}$. We can regard $X$ itself as an
element of the standard RBA and then we have
$$
(R X)^{[n]} = (0, M_{[n]}^1X, M_{[n]}^2X, \dots, M_{[n]}^lX,
\dots), \quad n\ge1,
$$
where $M_{[n]}^{l}$ is at the $(l+1)$-th position in the sequence.
Similarly
$$
(R X)^{\{n\}} = (0, M_{\omega_n}^1X, M_{\omega_n}^2X, \dots,
M_{\omega_n}^lX, \dots), \quad n\ge1.
$$

\begin{prop}
The elements $(RX)^{[n]}$ generate freely a subalgebra of~$\A$
(respectively generate freely a subalgebra of~$\A_R$).
\end{prop}

The proofs are omitted for the sake of brevity; the first uses the
observation that, for $l$ big enough, we find
$\Sup(M_{[n_1]}^l\cdots M_{[n_k]}^l)> \Sup(M_{[m_1]}^l\cdots
M_{[m_j]}^l)$ with $n_1+\cdots+n_k=m_1+\cdots+m_j$ iff the
sequence $(n_1,\ldots,n_k)$ is smaller than the sequence $(m_1,
\ldots,m_j)$ in the lexicographical ordering. The second is a bit
more involved.

\smallskip

The algebra~${\rm NCQSym}$ of quasi-symmetric functions in
noncommuting variables is naturally provided with a Hopf algebra
structure~\cite{bergeron05}.  On the elementary quasi-symmetric
functions $M_{[n]}$, the coproduct~$\Delta$ acts as on a sequence
of divided powers: $\Delta \bigl(M_{[n]}\bigr) = \sum_{i=0}^n
M_{[i]}\otimes M_{[n-i]}$.  Thus the $M_{[n]}$ generate a free
subalgebra of ${\rm NCQSym}$ naturally isomorphic as a Hopf
algebra to the classical descent algebra, which is a convolution
subalgebra of the endomorphism algebra
of~$T(X)$~\cite{reutenauer93} ---or equivalently, to the algebra
of noncommutative symmetric functions (NCSF)~\cite{gelfandetal95}.
The same construction goes over to the free algebras over the
$(RX)^{[n]}$ for the pointwise product and the Rota--Baxter double
product $\ast_R$.  The first one is naturally provided with a
cocommutative Hopf algebra structure for which the $(RX)^{[n]}$s
form a sequence of divided powers, that is:
$$
\Delta \bigl((RX)^{[n]}\bigr) = \sum\limits_{0\leq m\leq n}
(RX)^{[m]}\otimes (RX)^{[n-m]};
$$
this is just the structure inherited from the Hopf algebra structure
on ${\rm NCQSym}$.  We call this algebra the free noncommutative
Spitzer (Hopf) algebra on one generator, or the \textit{Spitzer
algebra} for short, and write $\mathcal S$ for it.  When dealing with
the $\ast_R$ product, the right subalgebra to consider, as it will
emerge soon, is the free algebra freely generated by the
$(RX)^{[n]}X$.  We also make it a Hopf algebra by requiring the free
generators to form a sequence of divided powers, that is
\begin{equation*}
\Delta_\ast\big((RX)^{[n]}X\big) = 1\otimes (RX)^{[n]}X +
\sum\limits_{0\leq m\leq n-1} (RX)^{[n-m-1]}X \otimes (RX)^{[m]}X +
(RX)^{[n]}X \otimes 1.
\end{equation*}
Thus it is convenient to set $(RX)^{[-1]}X=1$.  We call this Hopf
algebra the \textit{double Spitzer algebra}, and write $\mathcal C$
for it.  We shall need the antipode for both Hopf algebras.  For this,
recourse to Atkinson's theorem~\cite{atkinson1963} seems the simplest
method.  Recall that we assume $\theta=1$.

\begin{thm} \rm{({\bf{Atkinson}}~\cite{atkinson1963})}
\label{Atkinson1}
Let $(A,R)$ be a unital Rota--Baxter algebra.  Fix $a\in A$ and let
$x$ and~$y$ be defined by $x=\sum_{n\in\N}t^n (Ra)^{[n]}$ and
$y=\sum_{n\in\N}t^n(\tilde{R}a)^{\{n\}}$, that is, as the solutions of
the equations
\begin{equation*}
x = 1 + tR(x\,a) \sepword{\rm and} y = 1 + t\tilde R(a\,y),
\end{equation*}
in~$A[[t]]$. We have the following factorization
\begin{equation*}
x\big(1 + at\big)y = 1, \sepword{\rm so that} 1 + a t = x^{-1}y^{-1}.
\end{equation*}
\end{thm}

\begin{cor}
\label{inverseAtkinson}
Let $(A,R)$ be an associative unital Rota--Baxter algebra.  Fix $a\in
A$ and assume~$x$ and~$y$ to solve the equations in the foregoing
theorem.  The inverses~$x^{-1}$ and~$y^{-1}$ solve the equations
\begin{equation*}
x^{-1} = 1 - tR(a\,y) \sepword{\rm and} y^{-1} = 1 - t\tilde R(x\,a),
\end{equation*}
in~$A[[t]]$.
\end{cor}

One checks $xx^{-1}=x^{-1}x=1$ by using the definitions and the
Rota--Baxter property.  Similarly for~$y^{-1}$.

\begin{cor}
\label{cor:RBantipode}
The action of the antipode~$S$ on the Spitzer algebra~$\mathcal S$,
is given by
\begin{equation*}
S\bigl((RX)^{[n]}\bigr) = -R\bigl(X({\tilde R}X)^{\{n-1\}}\bigr).
\end{equation*}
\end{cor}

Indeed, the Spitzer bialgebra is naturally graded.  The series
$\sum_{n\in\N}(RX)^{[n]}$ is a group-like element in~$\mathcal S$.
The inverse series computes the action of the antipode on the
terms of the series.  The corollary follows, since
$$
\Big(\sum\limits_{n\in\N}(RX)^{[n]}\Big)^{-1} = 1 - R\Big(X
\big(\sum\limits_{n\in\N}({\tilde R}X)^{\{n\}}\big)\Big).
$$

\begin{cor}
\label{cor:RBdoubleAntipode} The action of the antipode~$S$ on the
double Spitzer algebra~$\mathcal C$ is given by
\begin{equation}
S\big((RX)^{[n]}X\big) = -\big(X({\tilde R}X)^{\{n\}}\big).
\label{eq:aleluya}
\end{equation}
\end{cor}

For the proof, one can observe that the operator $R$
induces an isomorphism of free graded algebras between~$\mathcal C$
and~$\mathcal S$ (which is the identity on scalars).  That is, for any
sequence of integers $i_1,\dots,i_k$, we have:
$$
R\bigl((RX)^{[i_1]}X \ast_R \cdots \ast_R (RX)^{[i_k]}X\bigr) =
(RX)^{[i_1+1]} \cdots (RX)^{[i_k+1]}.
$$
Hence, this implies~\eqref{eq:aleluya}.

\begin{cor}
The free $\ast_R$ subalgebras of $A$ generated by the
$(RX)^{[n]}X$ and the $X({\tilde R}X)^{\{n\}}$ are canonically
isomorphic.  The antipode exchanges the two families of
generators.  In particular, the $X({\tilde R}X)^{\{n\}}$ form also
a sequence of divided powers in the double Spitzer algebra.
\end{cor}


\section{Enter the Dynkin map}
\label{sect:nclogder}

The Dynkin operator is usually defined as the multilinear map from
an associative algebra $B$ into itself given by the left-to-right
iteration of the associated Lie bracket,
$$
D(x_1,\dots,x_n) = [\cdots[[x_1, x_2], x_3]\cdots\!, x_n],
$$
where $[x,y]:=xy-yx$.  Specializing to $B=T(X)$, the Dynkin operator
can be shown to become a quasi-idempotent ---that is, its action on an
homogeneous element of degree~$n$ satisfies $D^2=nD$.  The associated
projector $D/n$ sends $T_n(X)$ to the component of degree~$n$ of the
free Lie algebra over $X$, see the monograph~\cite{reutenauer93}.
Now, $D$ can be rewritten in purely Hopf algebraic terms as $S\star
N$, where $N$ is the grading operator and~$\star$ the convolution
product in $\End(T(X))$.  This definition generalizes to any graded
connected cocommutative or commutative Hopf algebra~\cite{patreu2002}.
One actually deals there with a more general phenomenon, namely the
possibility to define an action of the classical descent algebra on
any graded connected commutative or cocommutative Hopf
algebra~\cite{patras1994}.

\begin{thm}
\label{thm:invDynkinRB} Let $H$ be an arbitrary graded connected
cocommutative Hopf algebra over a field of characteristic zero.  The
Dynkin operator $D\equiv S\star N$ induces a bijection between the
group~$G(H)$ of group-like elements of $H$ and the Lie
algebra~$\Prim(H)$ of primitive elements in~$H$.  The inverse morphism
from~$\Prim(H)$ to~$G(H)$ is given by
\begin{equation}
\label{eq:invDynkin} h = \sum\limits_{n\in\N}h_n \longmapsto
\Gamma(h):=\sum\limits_{n\in\N} \sum\limits_{i_1+...+i_k = n,
\atop i_1,\dots,i_k>0}\,\frac{h_1 \cdots h_k}{i_1(i_1 + i_2)
\cdots (i_1 + \cdots + i_k)}.
\end{equation}
\end{thm}

This corresponds to Theorem~4.1 in our~\cite{eunomia2006},
establishing the same formula for characters and infinitesimal
characters of graded connected commutative Hopf algebras.  The proof
follows from the one in that reference by dualizing the notions and
identities, and can be omitted.  In the particular case where $H$ is a
free associative algebra over a set of graded generators
$y_1,\ldots,y_n,\ldots$ and $H$ is provided with the structure of a
cocommutative Hopf algebra by requiring the $y_i$ to be a sequence of
divided powers, the images of the generators $y_i$ under the action
of~$D$ forms a sequence of primitive elements of~$H$ that generate
freely $H$ as an associative algebra.  This result is a direct
consequence of our theorem.  Two particular examples of such a
situation are well known.  If $H$ is the NCSF Hopf algebra, then $H$
is generated as a free associative algebra by the complete homogeneous
NCSF, which form a sequence of divided powers, and the corresponding
primitive elements under the action of the Dynkin operator are known
as the power sums NCSF of the first kind~\cite{gelfandetal95}.
Second, in the classical descent algebra the abstract Dynkin operator
sends the identity of~$T(X)$ to the classical Dynkin operator.  This
was put to use in~\cite{reutenauer93} to rederive classical identities
of the Lie type.

We contend that the same machinery can be used to rederive the already
known formulas for commutative RBAs, and moreover prove new formulas
in the noncommutative framework.  We compute inductively the action
of~$D$ on the generators of~$\mathcal C$; that will give the action on
the generators of~$\mathcal S$, too.  Let us denote for the purpose
by~$\pi_*$ the product on~$\mathcal C$.  Using $N(1)=0$ and $N(X)=1$,
there follows $D((RX)^{[0]}X) = (S \star N)(X) = \pi_* \circ (S
\otimes N)\Delta_*(X)=\pi_* \circ (S \otimes N)(X \otimes 1 + 1
\otimes X) = X$.  We then find:
\begin{align*}
&D\bigl((RX)^{[n-1]}X\bigr) = (S\star N)\bigl((RX)^{[n-1]}X\bigr)
= \pi_*\circ(S \otimes N) \Bigl( \sum\limits_{0 \le p \le
n}(RX)^{[p-1]}X \otimes (RX)^{[n-p-1]}X\Bigr)
\\
&= \sum\limits_{0\le p\le n}S\bigl((RX)^{[p-1]}X\bigr) *_R
N\bigl((RX)^{[n-p-1]}X\bigr)
\\
&= \sum\limits_{0\le p\le n-1}S\bigl((RX)^{[p-1]}X\bigr) \ast_R
N\bigl((RX)^{[n-p-1]}\,\bigr)X + S\bigl((RX)^{[p-1]}X\bigr) \ast_R
(RX)^{[n-p-1]}X
\\
&= \sum\limits_{0\le p\le n-1}S\bigl((RX)^{[p-1]}X\bigr) \ast_R
N\bigl((RX)^{[n-p-1]}\,\bigr)X - S\bigl((RX)^{[n-1]}X\bigr)
\\
&= \sum\limits_{0\le p\leq n-1}R\Bigl(S\bigl((RX)^{[p-1]}X\bigr)
\ast_R N\bigl((RX)^{[n-p-2]}X\bigr)\Bigr)X
\\
&\qquad\qquad - \sum\limits_{1\le p\leq
n-1}S\bigl((RX)^{[p-1]}X\bigr){\tilde
R}\Bigl(R\bigl(N(R^{[n-p-2]}X)\bigr)X\Bigr) -
S\bigl((RX)^{[n-1]}X\bigr).
\end{align*}
In the fourth line we used vanishing of $(S\star\id)((RX)^{[n-1]}X)$,
then $a\ast_R(Rb\,c)= R(a\ast_Rb)c-a{\tilde R}(Rb\,c)$; the rest
should be clear.  After further simple manipulations,
using~\eqref{eq:aleluya} it comes
$$
D\bigl((RX)^{[n-1]}X\bigr) = R\bigl(D((RX)^{[n-2]}X)\bigr)X +
X{\tilde R}\bigl(D((RX)^{[n-2]}X)\bigr).
$$
The calculation suggests we introduce a new product.

\begin{de}
Let $(A,R)$ be an associative Rota--Baxter algebra. Introduce the
binary operation
\begin{equation}
a \bullet_R b := Ra\, b - bRa - ba = [Ra, b] - ba = Ra\, b +
b{\tilde R}a, \label{eq:ayvaa}
\end{equation}
and the elements $c^{(n)}(a_1,\dots,a_n) := \bigl(\cdots \big((a_1
\bullet_R a_2) \bullet_R a_3 \big) \cdots \bullet_R a_{n-1}\bigr)
\bullet_R a_n$, for $n>1$, and $c^{(1)}(a_1):=a_1$.
\end{de}

We define $c^{(n)}(a)$ as the $n$-times iterated product
$c^{(n)}(a,\dots,a) = \big(\cdots \big((a \bullet_R a) \bullet a \big)
\cdots \bullet_R a \big) \bullet_R a$.  All these parenthesis are
unavoidable, as the composition $\bullet_R$ is not associative, see
next section.  As well we define $C^{(n)}(a) :=
R\big(c^{(n)}(a)\big)$.  In conclusion, we have proved

\begin{thm}
The action of the Dynkin operator, $D$, on the generators $(RX)^{[n]}$
of the Spitzer algebra (respectively on the generators $(RX)^{[n]}X$
of the double Spitzer algebra) is given by
$$
D((RX)^{[n]}) = C^{(n)}(X), \sepword{\rm respectively by} D((RX)^{[n]}X) =
c^{(n)}(X).
$$
\end{thm}

This immediately implies

\begin{cor}
\label{cor:madre-del-cordero}
We have the following identity in the Spitzer algebra $\mathcal{S}$
\begin{equation}
(RX)^{[n]} = \sum\limits_{i_1+ \cdots +i_k = n, \atop i_1,\dots,i_k>0}
\frac{C^{(i_1)}(X) \cdots C^{(i_k)}(X)}{i_1(i_1 + i_2) \cdots (i_1 +
\cdots + i_k)}.
\label{eq:key1}
\end{equation}
\end{cor}

\begin{cor}
\label{cor:padre-del-cordero} We have the following identity in
the double Spitzer algebra $\mathcal{C}$
$$
(RX)^{[n-1]}X = \sum\limits_{i_1 + \cdots +i_k=n, \atop
i_1,\dots,i_k>0} \frac{c^{(i_1)}(X) \ast_R \cdots \ast_R
c^{(i_k)}(X)}{i_1(i_1 + i_2) \cdots (i_1 + \cdots +i_k)}.
$$
\end{cor}

The corollaries follow readily from our
Theorem~\ref{thm:invDynkinRB} by applying the inverse Dynkin
map~(\ref{eq:invDynkin}).


\section{The generalized Bohnenblust--Spitzer identities}
\label{sect:NCBS}

If $(A,R)$ is a \textit{commutative} Rota--Baxter algebra of weight
one with Rota--Baxter operator~$R$, then on $A[[t]]$ the following
identity by Spitzer holds~\cite{baxter1960,spitzer1956}:
\begin{equation}
\sum_{m \in \mathbb{N}} t^m (Ra)^{[m]} =
\exp\bigl(R\log(1+at)\bigr). \label{eq:collateral-damage}
\end{equation}
In the framework of the commutative standard RBA this becomes
Waring's formula relating elementary and power symmetric
functions~\cite[Chapter~4]{Sagan}.
From~\eqref{eq:collateral-damage} follows
$$
n!\,(Ra)^{[n]} = \sum_\sigma(-1)^{n-k(\sigma)}Ra^{|\tau_1|}\,
Ra^{|\tau_2|} \cdots Ra^{|\tau_{k(\sigma)}|}.
$$
Here the sum is over all permutations $\sigma$ of $[n]$ and
$\sigma=\tau_1\tau_2\cdots\tau_{k(\sigma)}$ is the
decomposition of~$\sigma$ into disjoint cycles~\cite{rotasmith72}.
We denote by $|\tau_i|$ the number of elements in $\tau_i$. By
polarization one obtains
\begin{equation*}
\sum_\sigma R\Bigl( R\bigl( \cdots (Ra_{\sigma(1)}) a_{\sigma(2)}
\cdots \bigr) a_{\sigma(n)}\Bigr) =
\sum_\sigma(-1)^{n-k(\sigma)}R\Big(\prod_{j_1\in\tau_1}a_{j_1}\Big)
\cdots
R\Big(\prod_{j_{k(\sigma)}\in\tau_{k(\sigma)}}a_{j_{k(\sigma)}}\Big).
\end{equation*}
This leads to the classical formula~\cite{rotasmith72}
\begin{equation}
\sum_{\sigma} R\Bigl( R\bigl( \cdots (Ra_{\sigma(1)})
a_{\sigma(2)} \cdots  \bigr) a_{\sigma(n)}\Bigr) = \sum_{\pi \in
\mathcal{P}_n} (-1)^{n -|\pi|} \prod_{\pi_i \in \pi}(m_i-1)! \
R\Bigl(\prod_{j \in \pi_i}a_j\Bigr). \label{eq:BohnenblustSp}
\end{equation}
Here $\pi$ now runs through all unordered set partitions
$\mathcal{P}_n$ of $[n]$; by~$|\pi|$ we denote the number of
blocks in~$\pi$; and $m_i:=|\pi_i|$ is the size of the particular
block~$\pi_i$. Those are often called Bohnenblust--Spitzer
formulas.  The generalization to \textit{noncommutative}
Bohnenblust--Spitzer formulas springs from
Corollaries~\ref{cor:madre-del-cordero},
respectively~\ref{cor:padre-del-cordero}.  Moreover, we arrive at
the following theorem.

\begin{thm}
\label{thm:factorization}
Let $(A,R)$ be an associative Rota--Baxter algebra.  For $a_i \in A$,
$i=1,\dots,n$, we have
 \allowdisplaybreaks{
\begin{eqnarray}
&&\sum_\sigma R\Bigl(R\bigl(\cdots
(Ra_{\sigma(1)}) a_{\sigma(2)} \cdots \bigr) a_{\sigma(n)}\Bigr) =
\sum_\sigma R\Bigl(a_{\sigma(1)} \diamond_1 a_{\sigma(2)}
\diamond_2 \dots \diamond_n a_{\sigma(n)}\Bigr), \sepword{\rm
where}
\label{eq:ayvahombre}\\
&&a_{\sigma(i)} \diamond_i a_{\sigma(i+1)} =
\begin{cases}a_{\sigma(i)} \ast_R a_{\sigma(i+1)}, & {\rm max} \;
(\sigma(j)| j\le i) < \sigma(i+1)
\\
a_{\sigma(i)} \bullet_R a_{\sigma(i+1)}, & \sepword{\rm
otherwise;}
\end{cases} \nonumber
\end{eqnarray}}
furthermore consecutive $\bullet_R$ products should be performed
from left to right, and always before the $\ast_R$ product.
\end{thm}

The reader might wish to perform a few checks here.  One readily finds
$$
R\bigl(Ra_1\,a_2\bigr) + R\bigl(Ra_2\,a_1 \bigr) = Ra_1 \, Ra_2 +
R(a_2 \bullet_R a_1) = R\bigl(a_1 \ast_R a_2 + a_2 \bullet_R
a_1\bigr) = R\bigl(a_2 \ast_R a_1 + a_1 \bullet_R
a_2\bigr).
$$
This is a fancy way to write the Bohnenblust--Spitzer identity in
terms of the non-associative Rota--Baxter product $\bullet_R$ and
the associative Rota--Baxter double product~$*_R$. To check by
direct calculation that \allowdisplaybreaks{
\begin{eqnarray*}
    \sum_{\sigma\in S_3}
        R\Bigl(R\bigl(Ra_{\sigma(1)}\, a_{\sigma(2)}\bigr) a_{\sigma(3)}\Bigr)
        &=&  R(a_1 \ast_R a_2 \ast_R a_3)
       + R\bigr(a_1 \ast_R ( a_3 \bullet_R  a_2)\bigr)
       + R\bigr(a_2 \ast_R (a_3 \bullet_R a_1)\bigr)
\nonumber \\[-0.4cm]
    &&
       + R\bigr((a_2 \bullet_R a_1) \ast_R a_3\bigr)
       + R\bigr((a_3 \bullet_R a_2) \bullet_R a_1 \bigr)
       + R\bigr((a_3 \bullet_R a_1) \bullet_R a_2\bigr)
\nonumber \\[0.1cm]
    &=&
         Ra_1\, Ra_2\, Ra_3
       + Ra_1\, R(a_3 \bullet_R a_2)
       + Ra_2\, R(a_3 \bullet_R a_1)
\nonumber \\
    &&
       + R(a_2 \bullet_R a_1)\, Ra_3
       + R\bigl((a_3 \bullet_R a_2) \bullet_R a_1\bigr)
       + R\bigl((a_3 \bullet_R a_1) \bullet_R a_2\bigr)
\end{eqnarray*}}
is already somewhat tedious.  We give a practical rule for the
decomposition in Theorem~\ref{thm:factorization}.  Given any
permutation~$\sigma$ of~$[n]$, we place a vertical bar to the left
of $\sigma_{i+1}$ iff it is bigger than all numbers to its left.
For instance, for $n=3$ we obtain in the one-line notation the
`cut permutations' $(1|2|3), (21|3), (312), (1|32), (321),
(2|31)$.  The cuts indicate where the $*_R$ products, if any,
should be placed.  Of course, as the left hand side
of~\eqref{eq:ayvahombre} is symmetrical in its arguments,
alternative rules could be devised.  For the decomposition of
$\sum_\sigma R\bigl(a_{\sigma(1)}R(a_{\sigma(2)} \cdots
Ra_{\sigma(n)})\cdots\bigr)$ our rule is: place a vertical bar to
the right of $\sigma_i$ iff it is smaller than all numbers to its
right.  For $n=3$ we then obtain the `cut permutations' $(1|2|3)$,
$(21|3)$, $(31|2)$, $(1|32)$, $(321)$, $(231)$; note the
differences. Moreover, in this case the $\bullet_R$ product is
defined by $aRb-Rb\,a-ba$ and consecutive $\bullet_R$ products are
performed from right to left.  As advertised, in the commutative
case, when $a\bullet_R b$ reduces to $-ab$, we recover the
classical Bohnenblust--Spitzer identities from any of the two
previous forms.


\section{Remarks and applications}
\label{sect:la gracia esta en los adverbios}

1.  Although the composition $\bullet_R$ in~(\ref{eq:ayvaa}) is
not associative, it is Vinberg or (left) pre-Lie. Recall that a
left pre-Lie algebra $V$ is a vector space, together with a
bilinear product $\bullet : V \otimes V \to V$, satisfying the
left pre-Lie relation
$$
(a \bullet b) \bullet c - a \bullet (b \bullet c) = (b \bullet a)
\bullet c - b \bullet (a \bullet c), \qquad a,b,c \in V.
$$
This is enough for the commutator $[a,b]:= a \bullet b - b \bullet
a$ to satisfy the Jacobi identity. Hence the algebra of
commutators $L_V$ is a Lie algebra, justifying the nomenclature.
Of course, every associative algebra is pre-Lie.
See~\cite{chapoliv2001} for more details on pre-Lie structures.

\begin{lem}
\label{lem:pre-LieRB} Let $(A,R)$ be an associative Rota--Baxter
algebra.  The binary composition~\eqref{eq:ayvaa} defines a left
pre-Lie structure on $A$, which we call left Rota--Baxter pre-Lie
product.
\end{lem}

The lemma follows by direct inspection.  It may also be related to
more recondite properties of RBAs~\cite{e2002}.  Let $(D,\ast)$ be
an associative algebra and assume that it is represented on
itself, from the left and from the right, with commuting actions.
We write $\succ$ and $\prec$ for the left and right actions,
respectively.  Assume moreover that we have $a \ast b = a \prec b
+ a \succ b$; then $D$ is by definition a dendriform dialgebra. In
detail, the dendriform properties are
\begin{equation}
(x \prec y )\prec z = x \prec (y \prec z + y \succ z); \; (x \succ
y )\prec z = x \succ (y \prec z); \; (x \prec y + x \succ y) \succ
z = x \succ (y \succ z). \label{eq:DD}
\end{equation}
Conversely, the latter relations are enough to ensure
associativity of $(D,\ast)$.  We refer to~\cite{loday2001} for
information on the subject.

Now, $D$ gives rise to a pre-Lie algebra and, in two different
ways, to the same Lie algebra. The pre-Lie algebra structure is
given by $x \bullet y := x \succ y - y \prec x$. As observed
already in~\cite{e2002}, generalizing an observation made by
Aguiar for the weight-zero case~\cite{aguiar2000}, the notion
applies in particular to weight $\theta \neq 0$ RBAs, since the
associative and pre-Lie products $\ast_R$ and $\bullet_R$,
respectively, are composed from sums and differences of the binary
operations
$$
a \prec_R b := - a \tilde{R}(b) \sepword{\rm and} a \succ_R b :=
R(a)b,
$$
that satisfy equations~\eqref{eq:DD} and define therefore a
dendriform dialgebra structure on any associative Rota--Baxter
algebra.  In the case of the Rota--Baxter pre-Lie composition, we
find
\begin{equation}
[a, b]_{\bullet_R} = [R(a), b] + [a, R(b)] + \theta [a, b] = [a,
b]_{\ast_R}. \label{eq:RBLiedouble}
\end{equation}

\begin{prop}
Let $(A,R)$ be an associative Rota--Baxter algebra.  The left
pre-Lie algebra $(A,\bullet_R)$ with the left Rota--Baxter pre-Lie
product is a Rota--Baxter pre-Lie algebra of the same weight, with
Rota--Baxter map~$R$.
\end{prop}

The proof of this is left as an exercise.

\smallskip

2.  It should be obvious now that, in the language of
NCSF~\cite{gelfandetal95}, if $X_a(t):=\sum_{n=0}^\infty
t^n(Ra)^{[n]}$ solves the initial value problem $d/dt\,X_a(t) =
X_a(t)\, \psi_a(t)$, $X_a(0)=1$, then $\psi_a(t) :=
\sum_{n>0}^\infty t^{n-1}C^{(n)}(a)$.

3.  The formulae developed in this paper actually apply without
restriction to any associative RBA, in particular to the solution
of differential equations.  We actually drew inspiration for this
paper from that subject: mainly from the path-breaking papers by
Lam~\cite{lamliu97,lam1998} and recent work by two of
us~\cite{Kalliope}.  To reestablish general weight in the pre-Lie
product formulas amounts simply to replace in~(\ref{eq:ayvaa}) the
product $ba$ by $\theta ba$, and thus the case $\theta=0$ is
included in our considerations.  In fact,
Corollary~\ref{cor:madre-del-cordero} yields the most efficient
way to organize the terms coming from the standard methods to
solve differential equations, the Dyson--Chen expansion and the
Magnus series.  Lam did obtain our formulas for~$(Ra)^{\{n\}}$ for
the case $\th=0$; part of the magic of the subject is how little
needs to be changed when $\theta \ne 0$.  It is worth mentioning
that this arose from the need to prove deep theorems with strong
physical roots, on approximations to quantum chromodynamics.  In
respect to the previous remark, if we define the Magnus series
coefficients $K_n$ by $d/dt\log X_a(t) = \sum_{n>0}^\infty t^n
K_n(a)$, then the relation between the $C^{(n)}$ and the $K_n$ is
precisely the relation between power sums NCSF of the first and of
the second kind~\cite{gelfandetal95}.  The advantage of writing
the Magnus series in this way has been recently recognized by the
practitioners~\cite{oteoros2000}.  Eventually, pointing to the
following remark we should underline that the NCSF picture implies
an exponential solution to Atkinson's recursion in
Theorem~\ref{Atkinson1}.

4.  It would be nice to be able to derive the new
Bohnenblust--Spitzer identities at one stroke from an equation
like the commutative Spitzer formula~\eqref{eq:collateral-damage}.
One of us participated in an attempt in this direction a few years
ago by~\cite{egk2005}, with the net result that in the
noncommutative case $\sum_m t^m (Ra)^{[m]}$ is still a functional
of $\log(1+at)$, through a non-linear recursion (for which
existence and unicity were proven) called, for want of a better
name, the Baker--Campbell--Hausdorff recursion, e.g.
see~\cite{egm2006}.  In practice, work with this functional was
painful.  There is a direct link between that recursion and the
Magnus expansion.  Explicit expressions for all the terms in the
latter are known; and so we are now forced to conclude that the
`solution' to the Baker--Campbell--Hausdorff recursion has been
staring at us for a while.  However, these formulas are rather
clumsy and will be presented elsewhere; the matter is under
investigation.

5. As shown in~\cite{eunomia2006}, the Dynkin operator is a key
ingredient for the mathematical understanding of the combinatorial
processes underlying the Bogoliubov recursion for renormalization
in perturbative quantum field theory. Use of general Spitzer-like
identities for noncommutative Rota--Baxter algebras is bound to
deepen this algebraic understanding of renormalization. From the
foregoing remarks it is clear that one can solve completely the
Bogoliubov recursion with this kind of Lie algebraic tools; this
will appear in a forthcoming work~\cite{Hit2007}.


\vspace{0.2cm}

\textbf{Acknowledgements}

The first named author acknowledges greatly the support by the
European Post-Doctoral Institute.  He also thanks Laboratoire
J.~A.~Dieudonn\'e at Universit\'e de Nice Sophia-Antipolis and the
Institut for Theoretical Physics at Bielefeld University for warm
hospitality.  JMG-B acknowledges partial support from CICyT,
Spain, through grant~FIS2005-02309.  The present work received
support from the ANR grant AHBE~05-42234.



\begin{thebibliography}{99}
\bibitem{bergeron05}
    N. Bergeron and M. Zabrocki,
    ``The Hopf algebras of symmetric functions and quasisymmetric
    functions in non-commutative variables are free and cofree'',
    \texttt{ArXiv:math.CO/0509265}.
\end{thebibliography}
\end{document}